\documentstyle[12pt]{article}
\if@twoside  m
\else
\fi
 1
 1
\font\sr = msbm10 scaled \magstep 1

\def\cross{{\triangleright\!\!\!<}}

\def\vaca(\langle\cdot\rangle\}
\def\diag{{\rm diag}\,}

\begin{document}
\begin{titlepage}
\begin{flushright}
CERN-TH/95-174\\
LBL-37416\\
UCB-PTH-95/18\\
q-alg/9508001\\
\end{flushright}
\vspace{.1cm}
\begin{center}
{\Large\bf Quantized Lax Equations and Their Solutions}\\
\vspace{1.3cm}
{ B. Jur\v co}\\
{\em CERN, Theory Division}\\
{\em CH-1211 Geneva 23, Switzerland}\\
\vspace{.5cm}
{M. Schlieker\footnote{This work was supported in part by the
Director, Office of Energy
Research, Office of High Energy and Nuclear Physics, Division of
High
Energy Physics of the U.S. Department of Energy under Contract
DE-AC03-76SF00098 and in part by the National Science Foundation
under grant PHY-90-21139.}$^,$\footnote{Supported in part by a
Feodor-Lynen Fellowship.}}\\
{\em Theoretical Physics Group}\\
{\em Lawrence Berkeley Laboratory}\\
{\em University of California}\\
{\em Berkeley, CA 94720, USA}

\end{center}
\vspace{1cm}
\begin{center}
{\bf Abstract}
\end{center}
{\small Integrable systems on quantum groups are investigated. The
Heisenberg equations possessing the Lax form are solved in terms of
the solution to the factorization problem on the corresponding
quantum group.}
\vspace{4cm}\\
CERN-TH/95-174\\
July 1995
\end{titlepage}

\newpage

\setcounter{footnote}{0}
\setcounter{page}{1}
{\thispagestyle{empty}
{\bf Disclaimer}\vspace{.5cm}\\
This document was prepared as an account of work sponsored by the
United
States Government.  Neither the United States Government nor any
agency
thereof, nor The Regents of the University of California, nor any of
their
employees, makes any warranty, express or implied, or assumes any
legal
liability or responsibility for the accuracy, completeness, or
usefulness
of any information, apparatus, product, or process disclosed, or
represents
that its use would not infringe privately owned rights.  Reference
herein
to any specific commercial products process, or service by its trade
name,
trademark, manufacturer, or otherwise, does not necessarily
constitute or
imply its endorsement, recommendation, or favoring by the United
States
Government or any agency thereof, or The Regents of the University  %
of
California.  The views and opinions of authors expressed herein do
not
necessarily state or reflect those of the United States Government
or any
agency thereof of The Regents of the University of California and
shall
not be used for advertising or product endorsement
purposes.\vspace{.5cm}\\

{\it Lawrence Berkeley Laboratory is an equal opportunity
employer.}\vspace{.5cm}\\

\newpage

\section{Introduction}
The discovery of the inverse scattering method \cite{Miura} was a
real breakthrough in theory of the classical completely integrable
Hamiltonian systems, which goes back to the classical papers of Euler,
Lagrange, Liouville, Jacobi and others. The systematic way to
construct and solve completely integrable Hamiltonian systems using
the theory of Lie groups and their representations originated in the
works of Kostant \cite{K}, Adler \cite{Ad} and Symes \cite{Symes}; it
was further  developed by many other authors.
The invention of Lie-Poisson groups by Drinfeld \cite{DrLP} made it
possible to develop the general concepts underlying the theory of
classical
integrable systems and the integration of the corresponding equations
of motions \cite{SemenovDressing}. We refer the reader to review
papers \cite{Re-Sem}, \cite{Tenlect} and to the books \cite{FT},
\cite{Newell} most related to our discussion. One of the main general
results in the theory is the construction  of integrable Hamiltonian
systems possessing a Lie (Lie-Poisson in general) group of symmetries
and the expression of their solutions in
terms of the solution to the factorization problem on this group
\cite{SemenovDressing}. The concrete classes of models differ by the
type of symmetry group and by the type of factorization.
Within this approach, the fundamental methods (inverse scattering
method, algebro-geometric methods of solution) and the fundamental
notions of the soliton theory, such as $\tau$-function \cite{JM} and
Baker-Akhieser function \cite{DKN}, \cite{SW}, found their unifying and
natural group-theoretical explanation.

The theory of integrable models of quantum mechanics and quantum field
theory also made a remarkable progress within the quantum version of
the inverse scattering method, which goes back to the seminal Bethe
ansatz for solving the Heisenberg
spin chain. We refer the reader, for a review of related topics, to the
books \cite{Gaudin}, \cite{Korepin} or to the papers \cite{F}, \cite{KS},
\cite{Thacker}. This
development made it possible to introduce quantum groups \cite{drinfeld},
\cite{J}
algebraic objects, playing in the quantum case a role analogous to that of
the Lie groups in the classical theory. However, we were still
missing (with the exception of the quantum integrable systems with
discrete time evolution \cite{Resh}) the quantum analogue of the
factorization theorem for the solution to the Heisenberg equations of
motion of a quantum integrable system. However, we have to mention the
remarkable paper \cite{Maillet} in this relation.

This paper is an attempt to formulate a quantum version of the
factorization theorem. As in the classical case, there is a simple direct
proof, and also a more conceptual proof which gives a generalization of the
classical
construction (based on symplectic reduction)
due to Semenov-Tian-Shansky \cite{SemenovDressing} and to which we devote the
main text. We hope that this construction has an interest of its own.
It could be, for example, interesting in relation with the Bethe ansatz
(Proposition 4), and it might be useful for a proper
formulation of the quantum version of the ${\tau}$-function. We prefer to
describe the direct proof in Appendix 1.

Section 2 contains the construction of a quantum dynamical system on
a dual quasi-triangular Hopf algebra
$F$, with the Hamiltonian $h$ taken as an arbitrary co-commutative
function on $F$, and gives its Lax pair formulation. Section 3
introduces a larger quantum dynamical system on the corresponding
Heisenberg double $D_H$ with a very simple time evolution, such that
the original quantum dynamical system (under some additional
assumptions) can be identified with a reduction of it. Section 4
contains our main result concerning the solution of the quantum Lax
equation of Section 2. Section 5 gives the formulation of our
results in a form suitable for integrable quantum chains or (after
performing the continuous limit) integrable quantum field theories.
Appendix 1 is devoted to a direct proof of our main theorem.
Finally in Appendix 2, we give a possible formulation of the factorization
problem in the case of factorizable Hopf algebras.

The paper is written for physicists. So, for example, we are working formally
with a notion of a dual Hopf algebra, which would need more detailed
specification in the infinite-dimensional case. Further, all algebraic tensor
products used in the paper would have to be properly completed in the
infinite-dimensional case.
Correspondingly we do
not discuss the precise sense in which the universal elements, such as
R-matrix, T-matrix,etc., exist. Apart from this, all constructions of the paper
are still valid. For these subtleties
and for more information about quantum groups we refer the reader to the
existing monographs on the subject (e.g. \cite{Pressleyatal}).

\section{Quantum Lax pairs}
The starting point of the following investigation is the
quasi-triangular Hopf algebra $U$ and its dual Hopf algebra $F = U^*$.

We shall use the standard notation: $m, \Delta$, $S$ and $\varepsilon$
for product, coproduct, antipode and co-unit, respectively, in both
$U$ and $F$,  and also
the notation $\sum x_{(1)}\otimes x_{(2)}$ for the result of
coproduct $\Delta$ applied to $x$, in $U$ or $F$.
We start from
the commutation relation \cite{FRT}
\begin{equation}
R_{12} T_{1} T_{2} = T_{2} T_{1} R_{12} \,,\label{br}
\end{equation}
where $R \in U \otimes U$ is the universal R-matrix and $T$ is the
universal element in $U \otimes F$ (sometimes called universal
T-matrix).

In the following we will always use the notation like
$$T_1 T_2 \equiv T_{13} T_{23} \,,$$
so that, for instance, (\ref{br}) means an equality in $U\otimes U\otimes F$
and the indices $1$ and $2$ refer to the different copies of $U$ in
this triple tensor product.

We want to study the following quantum dynamical system on the
quantum group $F$. The Hamiltonian $h$ is taken to be a co-commutative
element in $F$; it holds
\begin{equation}
\Delta (h) = \sigma \Delta (h) \,,
\end{equation}
where $\sigma$ is the flip operation.
The set of all such elements form a commutative subalgebra in $F$
\cite{drinfeld}.
The quantum dynamics is given by the following Heisenberg equations
of motion:
\begin{eqnarray}
i \dot{T_2}& = &[ h , T_2 ] \cr
&= &\langle h\otimes id , T_1 T_2 - T_2 T_1 \rangle \cr
&=& \langle h\otimes id , T_1 T_2 - T_2 T_1 -
{(R_{12}^{\pm})}^{-1}T_1T_2 + T_2 {(R_{12}^{\pm})}^{-1}T_1 \rangle\,,
\end{eqnarray}
where we used the commutation relations on $F$ with the universal
elements
$R^+ = R_{21}$ and $R^- = {R_{12}}^{-1}$ and the co-commutativity of
$h$, and where the time derivative applies
to the second tensor-factor of $T\in U\otimes F$ belonging to $F$.
Therefore we can state the following proposition, generalizing the discussion
of \cite{Maillet}:
\proclaim Proposition 1.
The Heisenberg equations can be written in the Lax form
\begin{eqnarray}
i \dot{T} &=& [ M^{\pm} , T ] \,, \hskip 0.5cm \mbox{where} \cr
M_2^{\mp} &=& \langle h \otimes id , ((1- {(R_{12}^{\pm})}^{-1})
T_1\rangle \,.\label{Lp}
\end{eqnarray}

In order to construct a solution for this set of equations we
shall consider in the next sections a quantized version of the
construction by
Semenov-Tian-Shansky in \cite{SemenovDressing}.

\section{Dynamics on the quantum Heisenberg double}

The quantum Heisenberg double $D_H$ (quantum cotangent bundle of $F$)
is a smash product algebra $D_H=U\cross F$ \cite{Sweedler} defined
with the help of the left action of $U$ on $F$.

We shall use the following description of $D_H$ \cite{Sem}, \cite{Zum}:
\proclaim Proposition 2.
The Heisenberg double $D_H$ is defined by
the following relations
\begin{eqnarray}
R_{12} T_1T_2 &=& T_2 T_1 R_{12} \,, \cr
R_{21} L_1^{\pm} L_2^{\pm} &=& L_2^{\pm} L_1^{\pm} R_{21} \,, \cr
R_{21} L_1^+ L_2^- &=& L_2^- L_1^+ R_{21} \,,\cr
L_1^{\pm} T_2 &=& T_2R_{12}^{\pm}L_1^{\pm} \,.
\end{eqnarray}

The universal T-matrix $T$ as well as universal L-matrices
$L^{\pm}=R^{\pm}$ are understood as elements of $U\otimes D_H$ in the
above equalities.

We shall introduce one more element of $U\otimes D_H$, denoted as $Y$
and defined as

$$Y=L^+(L^-)^{-1}.$$

Now let us consider the quantum dynamical system on the Heisenberg
double $D_H$
with the Hamiltonian ${\cal H}$ chosen to be a Casimir of $U \subset
D_H$ of the form \cite{GZB}
\begin{equation}
{\cal H}=\mbox{\mbox{Tr}}^v_1 ( Y_{12}^{-1} D_1)\,, \label{Ham}
\end{equation}
where $D\in U$ is defined, with the help of the universal R-matrix
$R=\sum R^{(1)}\otimes R^{(2)}$, as $D=\sum R^{(1)}S(R^{(2)})$ and the
superscript $v$ indicates the trace in the first factor of $U\otimes U
\subset U\otimes D_H$ evaluated
in an arbitrary representation $v$ of $U$.

The Heisenberg equations on $D_H$ take the form
\begin{eqnarray}
\dot Y &=& 0\,, \cr
i\dot T &=& T\xi^{\cal H}, \label{te}
\end{eqnarray}
with
\begin{equation}
\xi^{\cal H} =
\mbox{\mbox{Tr}}^v_1(R_{12}^{-1}Y_1^{-1}R_{21}^{-1}D_{1}-R_{12}^{-1}Y_
1^{-1}R_{12}D_{1})\,.
\end{equation}\label{xi}
Again the time derivative in (\ref{te}) applies to the second
tensor-factors of
$T$ and $Y$ belonging to $D_H$.
Since $\xi^{\cal H}\in U\otimes U \subset U\otimes D_H$ is
evidently time-independent, these Heisenberg equations are solved
trivially
\begin{eqnarray}
 Y(t)&=& Y(0)\,, \cr
T(t) &=& T(0)\mbox{exp}(-it\xi^{\cal H})\,. \label{ts}
\end{eqnarray}

In the following we will assume that $U$ itself is a quantum double
$D(U_-)$
of some Hopf algebra $U_-$.

Therefore we have
$$
U = U_-\otimes U_+\,,\hskip 1cm \mbox{with}\,U_+ = U_-^{* op \Delta}
$$
as a linear space and coalgebra. Similarly , we can write

$$
F = F_-\otimes F_+\,, \hskip 1cm \mbox{with}\,F_{\pm} = U^*_{\pm}
$$
as a linear space and an algebra.
Correspondingly we have $R\in U_-\otimes U_+$ and the universal
element $T$ factorizes in $U \otimes F $ as \cite{FRT}
\begin{eqnarray}
T &=& \Lambda Z \,, \cr
\mbox {with}
\,\Lambda &\in& U_- \otimes F_- \,,\cr
\mbox {and}\, Z &\in& U_+ \otimes F_+ \,.\label{dec}
\end{eqnarray}

The commutation relations of the elements $Y$ and $Z$ assumed as
elements in $U\otimes D_H$ play a crucial role in the following.
They are given by the following lemma.

\proclaim Lemma 1.
The elements $Y$ and $Z$ commute in the following way
\begin{eqnarray}
R_{21}Y_1R_{12}Y_2&=&Y_2R_{21}Y_1R_{12}\,, \cr
R_{12}Z_1Z_2&=&Z_2Z_1R_{12}\,, \cr
Z_1 Y_1 Z_2 &=& Z_2 Z_1 Y_1 R_{12} \,.\label{zy}
\end{eqnarray}

{\em Proof.}
Only the last assertion is non-trivial. We shall omit the details of
the proof of this relation, which follows immediately from the
discussion of \cite{Jurco-Schlieker} (all arguments given there we need are
valid also in
the general situation of the present paper), if we keep in mind the
difference in the decomposition of the universal
T-matrix used there and the decomposition (\ref{dec}). The resulting
difference is that the element $Q$ used in
\cite{Jurco-Schlieker} does not appear in the commutation relations
at all.

In order to make contact with the quantum dynamical system described
in Section 2, we  need the following proposition.

\proclaim Proposition 3.
There exists an embedding of $F\hookrightarrow D_H$, which is an
algebra homomorphism, given by
\begin{equation}
\tilde{T} = Z Y^{-1} Z^{-1} \,.\label{emb}
\end{equation}
This means that the relation
\begin{equation}
R_{12} \tilde{T_1} \tilde{T_2} = \tilde{T_2} \tilde{T_1} R_{12} \,
\end{equation}
holds in $U\otimes H_D$.
We shall use the symbol $\tilde F$ for the image of this embedding.

{\em Proof.}
The proof is straightforward using the commutation relations
(\ref{zy}).

This embedding of the original quantum group $F$ in the Heisenberg
double
will be used later on to project down a solution (\ref{ts}) to the
Heisenberg equations (\ref{te})
in $H_D$ to a solution of the Lax equation (\ref{Lp}) on the original
quantum phase space $F$. For doing this the following
identification of the
Hamiltonians of the corresponding systems is important.
Our starting Hamiltonian $h$ on the quantum group $F$ of Section
2 was supposed to be a co-commutative element in $F$. In the case when $F$
as its own left comodule decomposes to a direct sum of all its
irreducible comodules (a coarse form of the Peter-Weyl theorem)
the most general co-commutative element
$h$ is of the form

\begin{equation}
h=\mbox{\mbox{Tr}}^v T\,, \label{ham}
\end{equation}
where the trace in the first factor of $T\in U\otimes F$ is taken in
an appropriate representation $v$ of $U$. For simplicity, we shall
assume in the following our Hamiltonian $h$ to be exactly of
this type.

\proclaim Proposition 4.
For any representation $v$ of $U$ the equality
\begin{equation}
\mbox{\mbox{Tr}}^v_1 (Y_1^{-1} D_1) = \mbox{\mbox{Tr}}^v_1 ( Z_1
Y_1^{-1} Z_1^{-1}) = \mbox{\mbox{Tr}}^v_1 ( \tilde{T_1}) \,
\end{equation}
holds in $D_H$.
Roughly speaking the embedding (\ref{emb}) sends the trace of $\tilde
T$
in any representation of $U$ to the quantum trace of $Y$ in the same
representation, and we can identify the Hamiltonian $h$ with the reduction
to the $\tilde F$ of the Hamiltonian ${\cal H}$.

{\em Proof.} It holds in any representation $v$ of $U$ that
\begin{equation}
\mbox{\mbox{Tr}}_1 (\hat{R}_{12}^{-1} D_1 ) = 1 \,.\label{d-id}
\end{equation}
Here and in the rest of the proof we assume that both copies of $U$
to which indices
$1$ and $2$ refer are taken in the representation $v$.

{}From the third relation in (\ref{zy}) we get

\begin{equation}
Y_1^{-1} Z_1^{-1} Z_2 P_{12} = Z_2\hat{R}_{12}^{-1} Y_2^{-1}
Z_2^{-1}
\,,
\end{equation}
where $P_{12}$ is the permutation operator in the representation $v$.
Now taking the quantum trace of this equation and using (\ref{d-id})
we obtain

\begin{equation}
\mbox{\mbox{Tr}}^v_1 Y_1^{-1} Z_1^{-1} Z_2 P_{12} D_1 = Z_2 Y_2
Z_2^{-1} \label{*} \,.
\end{equation}
Taking now the usual trace in the second tensor-factor (and
renaming the tensor-factors) yields the desired identity
\begin{equation}
\mbox{\mbox{Tr}}^v_1 (Y_1^{-1} D_1 ) = \mbox{\mbox{Tr}}^v_1 (Z_1
Y_1^{-1} Z_1^{-1}) \,.
\label{trace}
\end{equation}

\section{Solution to the Lax equation}

Now we can return to our dynamical system, on $\tilde T\in U\otimes
D_H$ governed
by the Hamiltonian ${\cal H}$ of the form (\ref{Ham}), of Section 3. There, we
constructed the solution (\ref{ts}) to the equations
of motion
of this quantum dynamical system; we showed that
there exists an embedding of the original quantum group $F$ in the
Heisenberg double $D_H$, such that the Hamiltonians of the corresponding
systems coincide after this embedding.
In this chapter we are going to use it to obtain a
solution to the quantized Lax equations (\ref{Lp}) on the quantum group $F$.

Let us denote as $g(t)$ the following element

$$g(t)=Z(0)\mbox{exp}(-it\xi^{\cal H})Z(0)^{-1}\in U\otimes D_H.$$

Now, the time evolution on $D_H$ is an algebra homomorphism, and so the
decomposition of $T(t)$ in $U \otimes F(t)$ in the same form as in
(\ref{dec})
makes sense:
\begin{eqnarray}
T(t) &=& \Lambda(t) Z(t)\,,\cr
\mbox {with}\,
\Lambda &\in& U_- \otimes F_-(t) \,,\cr
\mbox{and}\, Z &\in& U_+ \otimes F_+(t) \,.\label{dec1}
\end{eqnarray}
So as a consequence of (\ref{ts}) the element $g(t)$ can be  expressed
as
\begin{eqnarray}
g(t) &=& \Lambda(0)^{-1} T(t) Z(0)^{-1} \,,\cr
&=& \Lambda(0)^{-1} \Lambda(t)Z(t) Z(0)^{-1} \,.\label{LZ}
\end{eqnarray}
This gives us a decomposition of $g(t) \in U\otimes
D_H$, with
\begin{eqnarray}
g(t)&=&g_-(t)g_+(t), \hskip 0.5cm g_{\pm}\in U_{\pm}\otimes D_H\,,\cr
g_-(t)&=&\Lambda(0)^{-1} \Lambda(t), \hskip 0.5cm g_+(t)=Z(t)
Z(0)^{-1}. \label{decc}
\end{eqnarray}

We will now show that $g(t)$ and its factors $g_{\pm}(t)$ are actually
elements of $U\otimes \tilde F
\subset U\otimes D_H$.

Let us define
\begin{equation}
M = R_- (h_{(1)}) \otimes h_{(2)} - R_+ (h_{(1)}) \otimes h_{(2)} \in U\otimes
\tilde F,
\end{equation}
with
\begin{eqnarray}
R_+(x)&=&\sum\langle x, R^{(1)}\rangle R^{(2)}\,, \cr
R_-(x)&=&\sum S(R^{(1)})\langle x, R^{(2)}\rangle\,, \label{maps}
\end{eqnarray}
and demonstrate that
\begin {equation}
g(t)=\mbox{exp}(-itM(0)) \in U\otimes
\tilde F\,.\label{g}
\end{equation}
That this is really true follows from the definition of $\xi^{\cal
H}$ (\ref{xi}), co-commutativity of $h$
and the following chain of identities:
\begin{eqnarray}
\mbox{\mbox{Tr}}_1^v(R_{12}^{-1}Y_1^{-1}(R_{12}^{\pm})^{-1}D_1)&=&Z_2^
{-1}Z_2\mbox{\mbox{Tr}}_1^v(R_{12}^{-1}Y_1^{-1}(R_{12}^{\pm})^{-1}D_1)
\,, \cr
&=&Z_2^{-1}\mbox{\mbox{Tr}}_1^v(Y_1^{-1}Z_1^{-1}Z_2Z_1(R_{12}^{\pm})^{
-1}D_1)\,,\cr
&=&Z_2^{-1}\mbox{\mbox{Tr}}_1^v(Y_1^{-1}Z_1^{-1}(R_{12}^{\pm})^{-1}Z_1
D_1)Z_2\,,\cr
&=&Z_2^{-1}\mbox{\mbox{Tr}}_1^v(Z_1Y_1^{-1}Z_1^{-1}(R_{12}^{\pm})^{-1}
)Z_2\,,\label{le}
\end{eqnarray}
where we used successively the third and the second relations of (\ref{zy}) and
the
relation
(\ref{*}).

Let us mention that $M=M^+-M^-$, with $M^{\pm}$ the elements of
$U\otimes \tilde F$
entering the Lax equation (\ref{Lp}).
{}From the equality $\xi^{\cal H}=Z^{-1}MZ$ that we just proved, and from the
time independence of $\xi^{\cal H}$, we have
$$ M(t) = g_+(t)M(0)g_+(t)^{-1}. $$

Writing now
\begin{equation}
g_+(t)=Z(t)Z(0)^{-1}
=\mbox{exp}(-it(1\otimes{\cal H}))Z(0)\mbox{exp}(it(1\otimes{\cal H}))
Z(0)^{-1}
\end{equation}
and using the last equality (with a $-$ sign) in (\ref{le}) and Proposition 4
we get immediately
\begin{equation}
g_+(t)=\mbox{exp}(-it(1\otimes h)\mbox{exp}(-it(M^+(0) -1\otimes
h)).\label{fga}
\end{equation}
For $g_-(t)$ we get similarly
\begin{equation}
g_-(t)=\mbox{exp}(-it(1\otimes h-M^-(0)))\mbox{exp}(it(1\otimes h)),\label{fgb}
\end{equation}
which follows from (\ref{fga}) and (\ref{commzero}) (in Appendix 1).
This shows that indeed $g_{\pm}\in U_{\pm}\otimes \tilde F$, as we claimed.
Moreover $g_{\pm}$ are the unique solutions of the equations
\begin{equation}
i\dot g_+=M^+g_+ \label{M+}.
\end{equation}
and
\begin{equation}
i\dot g_-=-g_-M^-,\label{M-}
\end{equation}
with initial condition $g_{\pm}(0)=1$.
Starting now from:
\begin{eqnarray}
\tilde T (t) &=& Z(t) Y^{-1}(0) Z(t)^{-1}\,,\cr
&=& Z(t)Z(0)^{-1}\tilde T(0)Z(0) Z(t)^{-1}\,,\label{tau}
\end{eqnarray}
we arrive at the main result of this paper:

\proclaim Theorem 1. Let $U$ be a quasi-triangular Hopf
algebra and let $F$ be its dual Hopf algebra. Let $g(t)$ be given by
(\ref{g}), with the Hamiltonian $h$, taken to be any co-commutative element of
$F$, and let $U_{\pm}$ denote the ranges of the mappings $R_{\pm}$
(\ref{maps}). Then $g(t)$ can be factorized:
\begin{equation}
g(t)=g_-(t)g_+(t)\,,\label{fac}
\end{equation}
$g_{\pm}(t)\in U_{\pm}\otimes F$ given by (\ref{fga}),  (\ref{fgb}). Moreover
$g_{\pm}(t)$ are the unique solutions of
equations (\ref{M+}), (\ref{M-}), with initial conditions $g_{\pm}(0)=1$.
The element $T(t) \in U\otimes F$, given by
\begin{equation}
 T(t)= g_+(t)T(0)g_{+}(t)^{-1} = g_-(t)^{-1}T(0)g_-(t) \,,\label{main}
\end{equation}
solves the quantum Lax equation (\ref{Lp}). In the case of factorizable $U$ we
can interpret (\ref{fac})
as a well-formulated factorization problem in $U\otimes F$ (see Appendix 2).

Although we proved here Theorem 1 only in the special case of $U$ being
a quantum double and the Hamiltonian $h$ being of the form (\ref{ham}), we
formulated it more generally. We shall give a simple direct proof of Theorem
1 in full generality in
Appendix 1.

The second equality in (\ref{main}) is due to fact that $g(t)$
commutes with $T(0)$ in $U\otimes F$, which is easily seen, e.g. from
(\ref{Lp}).

To specify completely our quantum dynamical system, we have to choose
a representation $\pi$ of the quantum group $F$. The algebra of
quantum observables will be the image $\pi(F)$ of $F$ in the chosen
representation.
The time evolution of an observable $\pi(a)$, $a \in F$, will then be
given by
$\pi(a)(t)=\pi((\langle a\otimes id), T(t)\rangle)$.

\section{Lax equations for quantum chains}

In this section we will discuss how the above result modifies in the
case of a quantum spin chain. The algebra of observables
$F^{\otimes N}$ for the chain consists of $N$ independent copies
$F^n,\,n=1,2,...,N$, of the dual Hopf algebra $F$ of a
quasi-triangular Hopf algebra $U$. We also ste $N+1\equiv 1$.
We will denote as $L^n \in U\otimes F^n$ the copy of the
universal T-matrix corresponding to the site $n$. We reserve the
character $T$ for the quantum monodromy matrix
$$T= L^1 ...L^N\,\in U\otimes F^{\otimes N}\,.$$
Then we have the following relations in $U\otimes F^{\otimes N}$
\begin{eqnarray}
R_{12} L_1^i L_2^i &=& L_2^i L_1^i R_{12} \,,\cr
L_1^i L_2^j &=& L_2^j L_1^i \,,\hskip 0.5cm i\neq j \,. \label{commi1}
\end{eqnarray}
The quantum monodromy matrix
satisfies
\begin{equation}
R_{12} T_1 T_2 = T_2 T_1 R_{12}\,
\end{equation}
and for the partial products
$$
{\psi}^n = L^1 ...L^{n-1} \,, \hskip 0.5cm {\psi}^1=1\,,
$$
we obtain
\begin{equation}
R_{12} {\psi}_1^n {\psi}_2^n = {\psi}_2^n {\psi}_1^n
R_{12}\,.\label{commi2} \end{equation}

We will choose our Hamiltonian $h\in F^{\otimes N}$ as any element of
$F^{\otimes N}$ of the form
\begin{equation}
h=\langle (H\otimes id),T \rangle\,, \label{chaham}
\end{equation}
with co-commutative $H\in F$. Again such elements form a commutative
subalgebra in $F^{\otimes N}$.

\proclaim Proposition 5.
The Lax equations for the site $n$ have the following
form:
\begin{eqnarray}
i\dot{L}^n &=& M^{\pm n} L^n - L^n M^{\pm (n+1)}, \hskip 0.5cm
\mbox{where} \cr
M_2^{\mp n} &=& \langle H \otimes id , ((1- (R_{12}^{\pm})^{-1})
(\psi_1^n)^{-1} T_1 \psi_1^n \rangle \,.  \label{cLp}
\end{eqnarray}

This can be easily shown using the commutation relations
(\ref{commi1}), (\ref{commi2}) and co-commutativity of $H$ in the
same way as in Proposition 1.

Lax pair of Proposition 5 formalizes concrete examples of Lax
pairs known
for particular integrable quantum chains or integrable field
theoretical models
\cite{Korepin1}, \cite{Sklyanin}, \cite{Kulish-Sklyanin},
\cite{Maillet}, \cite{Zhang}, \cite{Sogo}.

To avoid a cumbersome notation we introduce again a notation
similar to that in the previous section:
$$\hat M= R_- (H_{(1)}) \otimes H_{(2)} - R_+ (H_{(1)}) \otimes H_{(2)}\in
U\otimes F$$
and
$$ M=\langle \hat M\otimes id, id\otimes T\rangle \in U\otimes
F^{\otimes N}.$$

The folowing modification of Theorem 1 can be proved analogically
as in the previous sections. The twisted Heisenberg double of
ref. \cite{Sem} should be used for this.
However, there is also a direct proof using Theorem 1.

\proclaim Theorem 2.
Let $U$ be a quasi-triangular Hopf algbera, $F$ its dual Hopf algebra.
Let us assume a quantum chain system as described above with the
Hamiltonian $h$ given in (\ref{chaham}), where $H$ is any co-commutative
element
of $F$. Then the elements $g^n(t)\in U\otimes F^{\otimes N}$:
\begin{equation}
g^n(t) =({\psi}^n(0))^{-1} \mbox{exp}(-itM(0)){\psi}^n(0)\,,\label{gn}
\end{equation}
can be decomposed as
\begin{equation}
g^n(t)=g^n_-(t)g^n_+(t)\,,
\end{equation}
with $g^n_{\pm}(t)\in U_{\pm}\otimes F^{\otimes N}$ given by
\begin{eqnarray}
g^n_+(t)&=&\mbox{exp}(-it(1\otimes h)\mbox{exp}(-it(M^{+n}(0) -1\otimes
h))\,,\cr
g^n_-(t)&=&\mbox{exp}(-it(1\otimes h-M^{-n}(0)))\mbox{exp}(it(1\otimes
h)).\label{fgl}
\end{eqnarray}
Moreover $g^n_{\pm}$ are the unique solutions of equations (\ref{M+}),
(\ref{M-}) (all entries indexed by $n$), with initial conditions
$g^n_{\pm}(0)=1$. The elements
$L^n(t)\in U\otimes F^{\otimes N}$
\begin{equation}
L^n (t)=  g_+^n (t) L^n (0) (g_+^{n+1} (t))^{-1} = (g_-^n (t))^{-1}
L^n (0) g_-^{n+1} (t)
\end{equation}
solve the chain Lax equations (\ref{cLp}). In the case of factorizable $U$,
elements $g_{\pm}$ can be thought of
as a solution to the factorization problem for $g$ as formulated in Appendix 2.

{\em Proof.}
Following the same reasoning as led to Proposition 1, we can establish
that the Heisenberg equations of motion for entries of the quantum monodromy
matrices
$$T^n=(\psi^n)^{-1} T \psi^n=L^n...L^NL^1...L^{n-1},$$
for chains obtained from the original one by a shift $(1,...,N)\mapsto
(n,...,N,1,...,n-1)$, are precisely of the form (\ref{Lp}), with
$M^{\pm}=M^{\pm n}$. So the time evolution of the quantum monodromy matrix
$T^n$ is given by Theorem 1, with $g(t)=\mbox {exp}(-it(M^{+n}(0)-M^{-n}(0))$,
which means that all elements $\mbox {exp}(-it(M^{+n}(0)-M^{-n}(0))\in U\otimes
F^{\otimes N}$ can be decomposed as claimed.

It remais only to show that $$\mbox
{exp}(-it(M^{+n}(0)-M^{-n}(0))=({\psi}^n(0))^{-1}
\mbox{exp}(-itM(0)){\psi}^n(0).$$
This is, however, a consequence of the co-commutativity of $H$ and the
following
equality:
\begin{equation}
\psi_1^n(R_{12}^{\pm})^{-1}
(\psi_1^n)^{-1} T_1=(\psi_2^n)^{-1}(R_{12}^{\pm})^{-1}
 T_1\psi_2^n,
\end{equation}
which easily follows from (\ref{commi1}), (\ref{commi2}).

The rest is trivial.
\vskip 0.3cm
In this paper we did not mention the dressing symmetries of the quantum
integrable systems at all. However, dressing symmetries
can be introduced in a way completely analogous to the classical case (for the
classical case see \cite{SemenovDressing}). This aspect of the theory
of quantum integrable systems will be discussed elsewhere.
\vskip 1cm

\noindent{\bf Acknowledgements}

The authors would like to thank Bruno Zumino
for many valuable discussions. B.J. wishes also to acknowledge
discussions
with H. Grosse, P. Kulish and N. Reshetikhin. He would like
to thank Professors Grosse and Zumino for their kind hospitality at
ESI and LBL, respectively.
\vskip 1cm
\noindent
{\bf Appendix 1: Direct proof of Theorem 1}
\vskip 0.6cm
As in the classical case there is a simple direct proof of Theorem 1
and hence also of Theorem 2, which as it follows from the discussion of the
preceding section, is a simple consequence of the Theorem 1.

Let $U$ be any quasi-triangular Hopf algebra, $F$ its dual Hopf algebra and
$U_{\pm}$ the range of the maps $R_{\pm}$ (\ref{maps}).
First of all let us mention that $g_{\pm}\in U_{\pm}\otimes F$, given by
(\ref{fga}) and (\ref{fgb}):
\begin{eqnarray*}
g_+(t)&=&\mbox{exp}(-it(1\otimes h)\mbox{exp}(-it(M^+(0) -1\otimes h))\,,\cr
g_-(t)&=&\mbox{exp}(-it(1\otimes h-M^-(0)))\mbox{exp}(it(1\otimes h)).
\end{eqnarray*}
solve the equations (\ref{M+}), (\ref{M-}) with initial condition
$g_{\pm}(0)=1$, for any co-commutative Hamiltonian $h$. This is easily
checked by a direct computation.
As an immediate consequence we find that $T(t)$ given by (\ref{main})
solve the Lax equations (\ref{Lp}).

Now we shall show that the elements
$R_- (h_{(1)}) \otimes h_{(2)}\in U_-\otimes F$ and  $R_+ (h_{(1)}) \otimes
h_{(2)}\in U_+\otimes F$ commute.
Let us compute
\begin{eqnarray*}
R_{20}^{-1}T_0R_{12}T_1&=&
R_{20}^{-1}R_{12}T_0T_1=R_{20}^{-1}R_{12}R_{10}T_1T_0R_{10}^{-1}\cr
&=&R_{10}R_{12}R_{20}^{-1}T_1T_0R_{10}^{-1}=R_{10}R_{12}T_1R_{20}^{-1}T_0R_{10}^{-1}.
\end{eqnarray*}
Dualizing the first and last term in the above chain of equalities, which take
place in $U\otimes U\otimes U\otimes F$, in the components $0$ and $1$ with
$h\otimes h\in F\otimes F$, and using the co-commutativity of the Hamiltonian
$h$,
we have
\begin{equation}
[R_- (h_{(1)}) \otimes h_{(2)},R_+ (h_{(1)}) \otimes h_{(2)}]=0.
\label{commzero}
\end{equation}
This shows that
\begin{equation}
g_-(t)g_+(t)=\mbox{exp}(-it(M^+(0)-M^-(0))).
\end{equation}
So we have proved Theorem 1 directly.
\vskip 1cm
\noindent
{\bf Appendix 2: Factorization problem}
\vskip 0.6cm
Here we make an attempt to formulate a quantum analogue of the factorization
problem from the classical case \cite{SemenovDressing} in the case where $U$ is
a factorizable Hopf algebra \cite{S-R}. Similarly to
\cite{S-R}, we can give in this case an equivalent description of
the algebra structure of the tensor product $U\otimes F$. We shall omit
details.

The claim is that as a linear space $U\otimes F=F^{(-)}\otimes F^{(+)}$, where
$F^{(\pm)}$ are subalgebras of $U\otimes F$, both as algebras isomorphic to
$F$. They are embedded into $U\otimes F$ via the
the following algebra morphisms:
\begin{eqnarray*}
{\cal R}_{\pm}: F&\hookrightarrow& U\otimes F\,, \cr
x&\mapsto&R_{\pm}(x_{(1)}) \otimes x_{(2)}\,,
\end{eqnarray*}
where $R_{\pm}$ are given by (\ref{maps}).
This vector space isomorphism can be made into an algebra isomorphism if the
commutation relations between the elements of the two copies $F^{(\pm)}$ of $F$
are introduced through
\begin{equation}
(1\otimes x)(y\otimes 1)= \langle R_{21}^{-1},y_{(1)} \otimes x_{(1)}\rangle
 y_{(2)} \otimes x_{(2)}\langle R_{21},y_{(3)} \otimes x_{(3)}\rangle\,,
\end{equation}
so that, as an algebra, $U\otimes F$ is isomorphic to a bicrossproduct of two
copies of $F$.

This means that any element $\alpha\in U\otimes F$ can be expressed as
\begin{equation}
\alpha = \sum \alpha_{i-}\alpha_{i+},
\end{equation}
with all  $\alpha_{i-}$ lying in the range of the map ${\cal R}_-$ and all
$\alpha_{i+}$ lying in the range of the map ${\cal R}_+$, respectively. All
${\alpha_{i\pm}}$ are given unambiguously.
It may happen that some  particular $\alpha\in U\otimes F$, if expressed
in this way, is a simple product of two factors
\begin{equation}
\alpha = \alpha_-\alpha_+,
\end{equation}
$\alpha_+$ being the image under ${\cal R}_+$ of a (unique)
invertible element $x\in F$ and $\alpha_-$ being image under ${\cal R}_-$ of
the inverse $x^{-1}$ of the same element $x$.
If this is the case, we shall refer to the unique elements
\begin{eqnarray}
\tilde \alpha_-&=&\alpha_-(1\otimes x)\,,\cr
\tilde \alpha_+&=&(1\otimes x^{-1})\alpha_+
\end{eqnarray}
as to the solution of the factorization problem for $\alpha \in U\otimes F$.

Clearly the elements $g_{\pm}$ and $g_{\pm}^n$ from Theorems 1 and 2 are, in
the case
of the factorizable $U$, solutions to the factorization problem for $g$ and
$g^n$, respectively.

Finally we have to note that in concrete examples it is possible to give an
alternative characterization of the factorization of elements $g$ or $g^n(t)$.
We shall
discuss this very briefly for typical example when our starting Hopf algebra is
the quantum double of a Yangian $Y$ \cite{drinfeld}: $U=D(Y)$. Other cases
are similar. Let $T_{\lambda}$ be the automorphism of $Y$ of \cite{drinfeld}.
We shall use the same notation $T_{\lambda}, \lambda\in$ {\sr C} for its
extension (via duality) to the full double.
Then the decomposition of $(T_{\lambda}\otimes id)g^n(t)$,
$$(T_{\lambda}\otimes id)g^n(t)=(T_{\lambda}\otimes
id)g_-^n(t)(T_{\lambda}\otimes id)g_+^n(t),$$
is uniquely determined by the assumption that $(T_{\lambda}\otimes
id)g_{\pm}^n(t)$ are regular as functions of $\lambda$ in {\sr C}$P_1\backslash
\{\infty\}$
and {\sr C}$P_1\backslash\{0\}$, respectively, and $(T_{\infty}\otimes
id)g_+^n=1$.

\end{document}